\documentclass[11pt,a4paper]{article}

\usepackage[english]{babel}
\usepackage{amsmath}
\usepackage{amsthm}
\usepackage{amssymb}
\usepackage{amsfonts}
\usepackage{color}
\usepackage[breaklinks]{hyperref}
\usepackage[all]{xy}

\renewcommand{\phi}{\varphi}
\renewcommand{\epsilon}{\varepsilon}
\newcommand{\td}{\,\mathrm{d}}
\newcommand{\calT}{\mathcal{T}}
\newcommand{\calH}{\mathcal{H}}

\newcommand{\const}{\textup{const}}
\newcommand{\vol}{\textup{vol}}
\newcommand{\RR}{\mathbb{R}}

\theoremstyle{plain}
\newtheorem{theorem}{Theorem}[section]
\newtheorem*{theoremA}{Theorem A}
\newtheorem*{theoremB}{Theorem B}

\newtheorem*{fact}{Fact}

\newtheorem{lemma}[theorem]{Lemma}
\newtheorem{corollary}[theorem]{Corollary}
\newtheorem*{corollaryC}{Corollary C}

\theoremstyle{definition}

\theoremstyle{remark}
\newenvironment{pf}{\begin{proof}[\sc Proof]}{\end{proof}}

\numberwithin{equation}{section}

\title{An integral formula for $L^2$-eigenfunctions of a fourth order Bessel-type differential operator}

\author{Toshiyuki Kobayashi\footnote{Partially supported by Grant-in-Aid for Scientific Research (B) (18340037), Japan Society for the Promotion of Science, and the Alexander Humboldt Foundation.}\\\small{Graduate School of Mathematical Sciences}\vspace{-0.1cm}\\\small{The University of Tokyo}\vspace{-0.1cm}\\\small{3-8-1 Komaba, Meguro, Tokyo, 153-8914, Japan.}\vspace{-0.1cm}\\\small{\textit{E-mail address:} \texttt{toshi@ms.u-tokyo.ac.jp}}\vspace{-0.1cm}\\\\Jan M\"ollers\footnote{Partially supported by the International Research Training Group 1133 ``Geometry and Analysis of Symmetries'', and the GCOE program of the University of Tokyo.}\\\small{Institut f\"ur Mathematik}\vspace{-0.1cm}\\\small{Universit\"at Paderborn}\vspace{-0.1cm}\\\small{Warburger Str. 100}\vspace{-0.1cm}\\\small{33098 Paderborn, Germany}\vspace{-0.1cm}\\\small{\textit{E-mail address:} \texttt{moellers@math.uni-paderborn.de}}\vspace{-0.1cm}\\\ }

\makeindex

\begin{document}

\maketitle

\begin{abstract}
We find an explicit integral formula for the eigenfunctions of a fourth order differential operator against the kernel involving two Bessel functions. Our formula establishes the relation between $K$-types in two different realizations of the minimal representation of the indefinite orthogonal group, namely the $L^2$-model and the conformal model.\\


\end{abstract}


\section{Introduction and statement of the results}

Let $\theta=x\frac{\td}{\td x}$ be the one-dimensional Euler operator. We consider the following representation of the Bessel differential operator
\begin{align*}
 \frac{1}{x^2}Q_\nu(\theta) &= \frac{\td^2}{\td x^2}+\frac{\nu+1}{x}\frac{\td}{\td x}-1, & \nu\in\mathbb{C},
\end{align*}
where $Q_\nu$ is the quadratic transform of the Weyl algebra $\mathbb{C}[x,\frac{\td}{\td x}]$ defined by
\begin{align*}
 Q_\nu(P) &= P(P+\nu)-x^2, & \mbox{for }P\in\mathbb{C}\left[x,\frac{\td}{\td x}\right].
\end{align*}
Our object of study is the $L^2$-eigenfunctions of the fourth order differential operator
\begin{align*}
 D_{\mu,\nu} :={}& \frac{1}{x^2}Q_\nu(\theta+\mu)Q_\nu(\theta)\\
 ={}& \frac{1}{x^2}\left((\theta+\mu)(\theta+\mu+\nu)-x^2\right)\left(\theta(\theta+\nu)-x^2\right).
\end{align*}
Throughout this article we assume that the parameters $\mu$ and $\nu$ satisfy the following integrality condition:
\begin{align}
 \mu\geq\nu\geq-1\mbox{ are integers of the same parity, not both equal to $-1$.}\label{eq:IntCond}
\end{align}
We then have the following fact (see \cite[Theorem A]{HKMM09a}):

\begin{fact}
The differential operator $D_{\mu,\nu}$ extends to a self-adjoint operator on $L^2(\mathbb{R}_+,x^{\mu+\nu+1}\td x)$ with only discrete spectrum which is given by
\begin{align*}
 \lambda_j^{\mu,\nu} &:= 4j(j+\mu+1), & j=0,1,2,\ldots
\end{align*}
The corresponding $L^2$-eigenspaces are one-dimensional.
\end{fact}

For instance, it is easily seen that the normalized $K$-Bessel function $\widetilde{K}_{\frac{\nu}{2}}(z):=(\frac{z}{2})^{-\frac{\nu}{2}}K_{\frac{\nu}{2}}(z)$ is an $L^2$-eigenfunction of $D_{\mu,\nu}$ for the eigenvalue $\lambda_0^{\mu,\nu}=0$.

The purpose of this article is to establish the following integral formula for $L^2$-solutions of the differential equation
\begin{equation}
 D_{\mu,\nu}u = \lambda_j^{\mu,\nu}u.\label{eq:DiffEq}
\end{equation}

\begin{theoremA}
Assume \eqref{eq:IntCond} and let $u$ be an $L^2$-solution of the differential equation \eqref{eq:DiffEq}. Then there exists a constant $A_j^{\mu,\nu}(u)$ such that for $\cos\vartheta+\cos\phi>0$:
\begin{multline}
 \int_0^\infty{u(x)\widetilde{J}_{\frac{\mu}{2}}(ax)\widetilde{J}_{\frac{\nu}{2}}(bx)x^{\mu+\nu+1}\td x}\\
 = A_j^{\mu,\nu}(u)\left(\frac{\cos\vartheta+\cos\phi}{2}\right)^{\frac{\mu+\nu+2}{2}}\widetilde{C}_j^{\frac{\mu+1}{2}}(\cos\vartheta)\widetilde{C}_{j+\frac{\mu-\nu}{2}}^{\frac{\nu+1}{2}}(\cos\phi),\label{eq:MainIntFormula}
\end{multline}
where we set $a:=\frac{\sin\vartheta}{\cos\vartheta+\cos\phi}$ and $b:=\frac{\sin\phi}{\cos\vartheta+\cos\phi}$. 
\end{theoremA}

Here $\widetilde{J}_\alpha(x)=\left(\frac{x}{2}\right)^{-\alpha}J_\alpha(x)$ denotes the normalized $J$-Bessel function and $\widetilde{C}_n^\lambda(x)=\Gamma(\lambda)C_n^\lambda(x)$ is the normalized Gegenbauer polynomial.

The differential equation \eqref{eq:DiffEq} has a regular singularity at $x=0$ with characteristic exponents $0$, $-\nu$, $-\mu$ and $-\mu-\nu$. Accordingly, the asymptotic behaviour of a non-zero $L^2$-solution $u$ of \eqref{eq:DiffEq} as $x\rightarrow0$ is of the following form (see \cite[Theorem 4.2 (1)]{HKMM09a}):
\begin{align}
 u(x) &\sim B_j^{\mu,\nu}(u)\times\left\{\begin{array}{ll}x^{-\nu}+o(x^{-\nu}) & \mbox{for }\nu>0,\\\log(\frac{x}{2})+o(\log(\frac{x}{2})) & \mbox{for }\nu=0,\\1+o(1) & \mbox{for }\nu=-1,\end{array}\right.\label{eq:Asymptotics}
\end{align}
with some non-zero constant $B_j^{\mu,\nu}(u)$. The constant $A_j^{\mu,\nu}(u)$ in Theorem A is determined by $B_j^{\mu,\nu}(u)$ as follows:

\begin{theoremB}
For any solution $u$ of \eqref{eq:DiffEq}:
\begin{align*}
 \frac{A_j^{\mu,\nu}(u)}{B_j^{\mu,\nu}(u)} &= (-1)^j\frac{j!2^{2\mu+\nu}\Gamma(\frac{\mu+2}{2})\Gamma(\frac{\mu-|\nu|+2}{2})\Gamma(j+\frac{\mu-\nu+2}{2})}{\Gamma(j+\frac{\mu-|\nu|+2}{2})\pi\Gamma(j+\mu+1)}\times \left\{\begin{array}{ll}\displaystyle\frac{2}{\Gamma\left(\frac{\nu}{2}\right)}&\mbox{for $\nu>0$,}\\\displaystyle-1&\mbox{for $\nu=0$,}\\\displaystyle-\frac{2}{\Gamma\left(\frac{\nu}{2}\right)}&\mbox{for $\nu=-1$.}\end{array}\right.
\end{align*}
\end{theoremB}

The proofs of Theorems A and B will be given in Sections \ref{sec:MinRep} and \ref{sec:ConstComp}, respectively.

In Section \ref{sec:Applications} we give some applications and discuss special values of Theorem A. One particularly interesting situation arises when both $\mu$ and $\nu$ are odd integers. In this case the solutions $u$ of \eqref{eq:DiffEq} can be expressed as
\begin{align*}
 u(x) = \const\times\left\{\begin{array}{ll}x^{-\nu}e^{-x}M_j^{\mu,\nu}(2x)&\mbox{for }\nu\geq1,\\e^{-x}M_j^{\mu,\nu}(2x)&\mbox{for }\nu=-1,\end{array}\right.
\end{align*}
for some polynomial $M_j^{\mu,\nu}$ (see \cite{HKMM09b}). For $\nu=\pm1$ these polynomials reduce to the classical Laguerre polynomials $M_j^{\mu,\pm1}(x)=L_j^\mu(x)$. Hence, for $\nu=\pm1$ the integral formula in Theorem A collapses to integral formulas for the Laguerre polynomials. Even these we could not trace in the literature.

\begin{corollaryC}
Let $\mu\geq1$ be an odd integer and $\cos\vartheta+\cos\phi>0$. Set $a:=\frac{\sin\vartheta}{\cos\vartheta+\cos\phi}$ and $b:=\frac{\sin\phi}{\cos\vartheta+\cos\phi}$. Then we have
\begin{multline*}
 \int_0^\infty{L_j^\mu(2x)\widetilde{J}_{\frac{\mu}{2}}(ax)\cos(bx)x^\mu e^{-x}\td x}\\
 = (-1)^j\frac{2^\mu}{\sqrt{\pi}}\left(\frac{\cos\vartheta+\cos\phi}{2}\right)^{\frac{\mu+1}{2}}\cos\left(j+\frac{\mu+1}{2}\right)\phi\ \widetilde{C}_j^{\frac{\mu+1}{2}}(\cos\vartheta)
\end{multline*}
and
\begin{multline*}
 \int_0^\infty{L_j^\mu(2x)\widetilde{J}_{\frac{\mu}{2}}(ax)\sin(bx)x^\mu e^{-x}\td x}\\
 = (-1)^j\frac{2^\mu}{\sqrt{\pi}}\left(\frac{\cos\vartheta+\cos\phi}{2}\right)^{\frac{\mu+1}{2}}\sin\left(j+\frac{\mu+1}{2}\right)\phi\ \widetilde{C}_j^{\frac{\mu+1}{2}}(\cos\vartheta).
\end{multline*}
\end{corollaryC}

Note that there appear $5$ parameters in the integral formula in Theorem A, namely $\mu$, $\nu$, $\vartheta$, $\phi$ and $j$. Our scheme (specialization of parameters, relation to representation theory) is summarized in the following diagram:
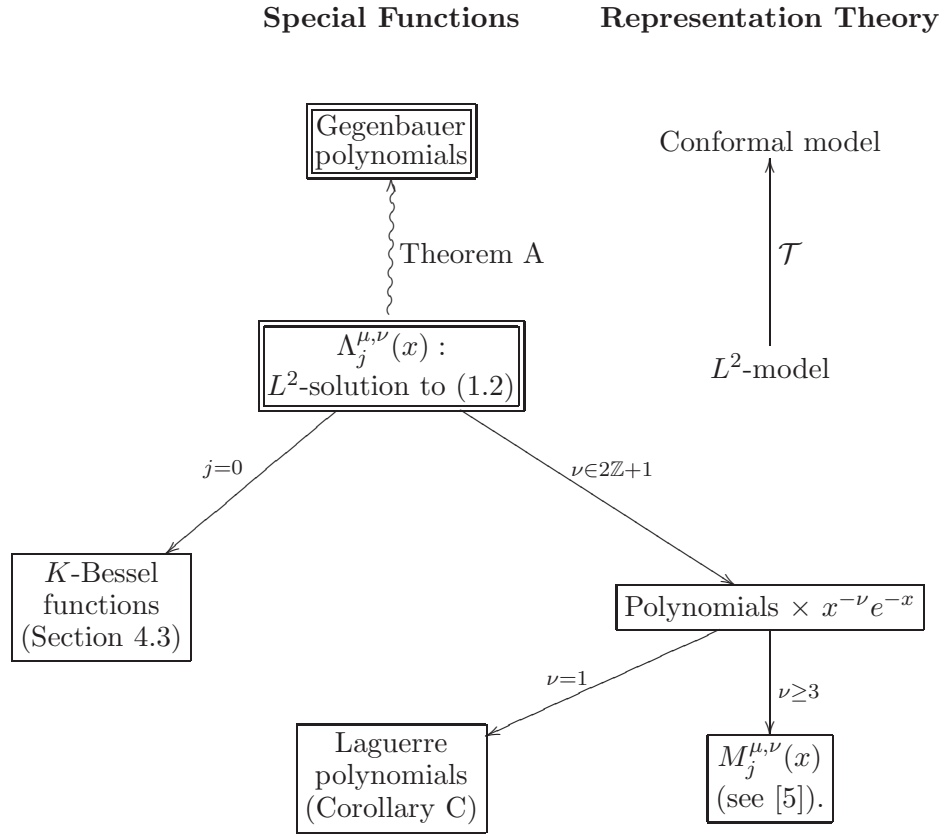
\[\xymatrix{
 & \mbox{\textbf{Special Functions}} & \mbox{\textbf{Representation Theory}}\\
 & *+[F=]{\stackrel{\mbox{Gegenbauer}}{\mbox{polynomials}}} & \mbox{Conformal model}\\\\
 & *+[F=]{\stackrel{\mbox{$\Lambda_j^{\mu,\nu}(x):$}}{\mbox{$L^2$-solution to \eqref{eq:DiffEq}}}} \ar@{~>}[uu]_{\mbox{Theorem A}} \ar[ddl]_{j=0} \ar[ddr]^{\nu\in2\mathbb{Z}+1} & \mbox{$L^2$-model} \ar[uu]_{\mbox{$\mathcal{T}$}}\\\\
 *+[F]\txt{\mbox{$K$-Bessel}\\\mbox{functions}\\\mbox{(Section \ref{sec:BottomCase})}} & & *+[F]{\mbox{Polynomials $\times\ x^{-\nu}e^{-x}$}} \ar[dl]_{\nu=1} \ar[d]^{\nu\geq3}\\
 & *+[F]\txt{\mbox{Laguerre} \\ \mbox{polynomials} \\ \mbox{(Corollary C)}} & *+[F]\txt{\mbox{$M_j^{\mu,\nu}(x)$} \\ \mbox{(see \cite{HKMM09b}).}}
}\]

Notation: $\mathbb{N}_0=\{0,1,2,\ldots\}$, $\mathbb{R}_+=\{x\in\mathbb{R}:x>0\}$.

\section{Two models for the minimal representation of the indefinite orthogonal group}\label{sec:MinRep}

In the proof of Theorem A we will use representation theory, namely two different models for the minimal representation of the indefinite orthogonal group $G=O(p,q)$ where $p\geq q\geq2$ and $p+q\geq6$ is even. These two models were constructed by T. Kobayashi and B. {\O}rsted \cite{KO03a,KO03c} and investigated further by T. Kobayashi and G. Mano \cite{KM07b}. This unitary representation is irreducible and attains the minimum Gelfand--Kirillov dimension among all irreducible unitary representations of $G$. In physics the minimal representation of $O(4,2)$ appears as the bound states of the Hydrogen atom, and incidentally as the quantum Kepler problem.

\subsection{The conformal model}
We begin with a quick review of the conformal model for the minimal representation of $G=O(p,q)$ from \cite{KO03a}.

We equip $M:=S^{p-1}\times S^{q-1}$ with the standard indefinite Riemannian metric of signature $(p-1,q-1)$ by letting the second factor be negative definite. Then $G$ acts on $M$ by conformal transformations. The solution space
\begin{align*}
 \mathcal{S}ol(\widetilde{\Delta}_M) := \left\{f\in C^\infty(M):\widetilde{\Delta}_Mf=0\right\}
\end{align*}
of the Yamabe operator $\widetilde{\Delta}_M=\Delta_{S^{p-1}}-\Delta_{S^{q-1}}-\left(\frac{p-2}{2}\right)^2+\left(\frac{q-2}{2}\right)^2$ is infinite-dimensional. Further, it is invariant under the \lq twisted action\rq\ $\varpi$ of $G$ and hence defines a representation. The minimal representation of $G$ is realized on the Hilbert completion
\begin{align*}
 \mathcal{H} := \overline{\mathcal{S}ol(\widetilde{\Delta}_M)}
\end{align*}
of $\mathcal{S}ol(\widetilde{\Delta}_M)$ with respect to a certain $G$-invariant inner product.

The maximal compact subgroup $K=O(p)\times O(q)$ of $G$ acts on $M$ as isometries, and the restriction of $\varpi$ to $K$ is given just by rotations. To see the $K$-types we recall the space of spherical harmonics
\begin{align*}
 \mathcal{H}^k(\mathbb{R}^n) &:= \left\{\phi\in C^\infty(S^{n-1}):\Delta_{S^{n-1}}\phi=-k(k+n-2)\phi\right\},
\end{align*}
or equivalently, the space of restrictions of harmonic homogeneous polynomials on $\mathbb{R}^n$ of degree $k$ to the sphere $S^{n-1}$. The orthogonal group $O(n)$ acts irreducibly on $\calH^k(\RR^n)$ for any $k$ by rotations in the argument. Then clearly
\begin{align*}
 \mathcal{H}^j(\mathbb{R}^p)\otimes\mathcal{H}^k(\mathbb{R}^q) \subseteq \mathcal{S}ol(\widetilde{\Delta}_M) \mbox{ if and only if } k=j+\frac{p-q}{2},
\end{align*}
and we put
\begin{align*}
 V^j &:= \mathcal{H}^j(\mathbb{R}^p)\otimes\mathcal{H}^{j+\frac{p-q}{2}}(\mathbb{R}^q), & j=0,1,2,\ldots,
\end{align*}
on which $K$ acts irreducibly. The multiplicity-free sum $\oplus_{j=0}^\infty{V^j}$ of irreducible representations of $K$ is dense in the Hilbert space $\mathcal{H}$.

Let $K'$ be the isotropy group of $K$ at $((1,0,\ldots,0),(0,\ldots,0,1))\in S^{p-1}\times S^{q-1}$. Then $K'\cong O(p-1)\times O(q-1)$. We write $\mathcal{H}^{K'}$ for the space of $K'$-fixed vectors.

\begin{lemma}
In each $K$-type $V^j$ ($j\in\mathbb{N}_0$) the subspace $V^j\cap\mathcal{H}^{K'}$ is one-dimensional and spanned by the functions
\begin{align}
 \psi_j:S^{p-1}\times S^{q-1}\rightarrow\mathbb{C}, (v_0,v',v'',v_{p+q-1})\mapsto\widetilde{C}_j^{\frac{p-2}{2}}(v_0)\widetilde{C}_{j+\frac{p-q}{2}}^{\frac{q-2}{2}}(v_{p+q-1}).\label{eq:KtypesConformal}
\end{align}
\end{lemma}

\begin{pf}
It is well-known that any $O(n-1)$-invariant spherical harmonic is a scalar multiple of the Gegenbauer polynomial
\begin{align*}
 S^{n-1}\ni(x_1,x')\mapsto\widetilde{C}_k^{\frac{n-2}{2}}(x_1),
\end{align*}
which shows the claim.
\end{pf}

\subsection{The $L^2$-model}

We recall the $L^2$-model (Schr\"odinger model) of the minimal representation of $G$ which is unitarily equivalent to $\varpi$ (see \cite{KM07b,KO03c}). Consider the isotropic cone
\begin{align*}
 C=\{(x',x'')\in\RR^{p-1}\times\RR^{q-1}:|x'|=|x''|\neq0\}\subseteq\RR^{p+q-2}.
\end{align*}
Then the group $G$ acts unitarily in a non-trivial way on the Hilbert space $L^2(C,\td\mu)$ and defines a minimal representation of $G$. Here $\td\mu$ is the $O(p-1,q-1)$-invariant measure on $C$ which is in bipolar coordinates
\begin{align*}
 \RR_+\times S^{p-2}\times S^{q-2}\stackrel{\sim}{\longrightarrow}C,\ (r,\omega,\eta)\mapsto(r\omega,r\eta)
\end{align*}
normalized by $\td\mu=\frac{1}{2}r^{p+q-5}\td r\td\omega\td\eta$. ($\td\omega$ and $\td\eta$ denote the Euclidean measures on $S^{p-2}$ and $S^{q-2}$, respectively.) The representation of the whole group $G$ on $L^2(C)$ does not come from the geometry $C$, but the action of the subgroup $K'$ is given by rotation in the argument. Hence the $K'$-invariant functions only depend on the radial parameter $r\in\RR_+$ and the space of $K'$-invariants in $L^2(C)$ is identified as $L^2(C)^{K'}\cong L^2(\RR_+,\frac{1}{2}r^{p+q-5}\td r)$.

Let $W^j$ be the $V^j$-isotypic component in $L^2(C)$. In this model it is more difficult to find explicit $K$-finite vectors. By highlighting $K'$-fixed vectors, the following result was proved in \cite[Section 8]{HKMM09a}:

\begin{lemma}
In each $K$-type $W^j$ ($j\in\mathbb{N}_0$) the subspace $W^j\cap L^2(C)^{K'}$ is one-dimensional and given by the radial functions
\begin{align}
 u(2r),\label{eq:KtypesL2}
\end{align}
where $u$ is an $L^2$-solution of \eqref{eq:DiffEq} with $\mu=p-3$, $\nu=q-3$.
\end{lemma}

\subsection{The $G$-intertwiner}

Let $\calT:L^2(C)\stackrel{\sim}{\longrightarrow}\mathcal{H}$ be the intertwining operator as given in \cite[Section 2.2]{KM07b}. It is the composition of the Fourier transform $\mathcal{S}'(\mathbb{R}^{p+q-2})\rightarrow\mathcal{S}'(\mathbb{R}^{p+q-2})$ and an operator coming from the conformal transformation from the flat indefinite Euclidean space $\mathbb{R}^{p-1,q-1}$ to $M$. For radial functions $f\in L^2(\RR_+,\frac{1}{2}r^{p+q-5}\td r)\cong L^2(C)^{K'}$ this operator can be written by means of the Hankel transform (cf. \cite[Lemma 3.3.1]{KM07b}):
\begin{multline}
 \calT f(v_0,v',v'',v_{p+q-1}) = \frac{1}{(v_0+v_{p+q-1})^{\frac{p+q-4}{2}}}\int_0^\infty{f(r)}\\
 \times{\widetilde{J}_{\frac{p-3}{2}}\left(\frac{2|v'|r}{v_0+v_{p+q-1}}\right)\widetilde{J}_{\frac{q-3}{2}}\left(\frac{2|v''|r}{v_0+v_{p+q-1}}\right)r^{p+q-5}\td r}\label{eq:RadialIntertwiner}
\end{multline}
for $(v_0,v',v'',v_{p+q-1})\in M$ with $v_0+v_{p+q-1}>0$. Since $\calT$ intertwines the actions of $G$ on both models, it clearly maps $K'$-invariant functions to $K'$-invariant functions and also preserves $K$-types, i.e. $\calT(W^j)=V^j$. Hence we obtain the following diagram:
\begin{align*}
 \xymatrix@R=1pc@C=0.00000001pc{ **[l] W^j \ar[r]^{\sim} \ar@{}[d]|{\displaystyle\cap\ \ \ } & **[r] V^j \ar@{}[d]|{\displaystyle\ \cap}\\ **[l] \calT:\ L^2(C) \ar[r]^{\sim} & **[r] \mathcal{H}^{K'}\\ **[l] L^2(C)^{K'} \ar[r]^{\sim} \ar@{}[u]|{\displaystyle\cup\ \ \ } & **[r] \mathcal{H}^{K'}. \ar@{}[u]|{\displaystyle\ \cup} }
\end{align*}

Now $W^j\cap L^2(C)^{K'}$ and $V^j\cap\mathcal{H}^{K'}$ are one-dimensional and we have formulas \eqref{eq:KtypesL2} and \eqref{eq:KtypesConformal} for their generators. Hence the intertwiner $\calT$ has to map the functions \eqref{eq:KtypesL2} to multiples of the functions \eqref{eq:KtypesConformal}. Thus, for any $L^2$-solution $u$ of \eqref{eq:DiffEq} and $v_0+v_{p+q-1}>0$ we obtain
\begin{multline}
 \int_0^\infty{u(2r)\widetilde{J}_{\frac{p-3}{2}}\left(\frac{2|v'|r}{v_0+v_{p+q-1}}\right)\widetilde{J}_{\frac{q-3}{2}}\left(\frac{2|v''|r}{v_0+v_{p+q-1}}\right)r^{p+q-5}\td r}\\
 =\const\cdot(v_0+v_{p+q-1})^{\frac{p+q-4}{2}}\widetilde{C}_j^{\frac{p-2}{2}}(v_0)\widetilde{C}_{j+\frac{p-q}{2}}^{\frac{q-2}{2}}(v_{p+q-1})\label{eq:IntFormulaConstant}
\end{multline}
Substituting $x=2r$, $\mu=p-3$ and $\nu=q-3$ and putting
\begin{align*}
 \cos\vartheta &= v_0, & \cos\phi &=v_{p+q-1},\\
 \sin\vartheta &= |v'|, & \sin\phi &=|v''|.
\end{align*}
we get \eqref{eq:MainIntFormula} with a certain constant $A_j^{\mu,\nu}$. This finishes the proof of Theorem A.

\section{A closed formula for the constants}\label{sec:ConstComp}

In this section we find an explicit constant for the integral formula in Theorem A, namely we give a proof of Theorem B. Our method uses the generating function of $L^2$-eigenfunctions of $D_{\mu,\nu}$.

Remember that we assume the integrality condition \eqref{eq:IntCond}. Let
\begin{align}
 G^{\mu,\nu}(t,x) &= \frac{1}{(1-t)^{\frac{\mu+\nu+2}{2}}}\widetilde{I}_{\frac{\mu}{2}}\left(\frac{tx}{1-t}\right)\widetilde{K}_{\frac{\nu}{2}}\left(\frac{x}{1-t}\right),\label{eq:GenFct}
\end{align}
where $\widetilde{I}_\alpha(z):=(\frac{z}{2})^{-\alpha}I_\alpha(z)$ and $\widetilde{K}_\alpha(z):=(\frac{z}{2})^{-\alpha}K_\alpha(z)$ denote the normalized $I$- and $K$-Bessel functions. Further, let $(\Lambda_j^{\mu,\nu}(x))_{j=0,1,2,\ldots}$ be the family of functions on $\mathbb{R}_+$ which has $G^{\mu,\nu}(t,x)$ as its generating function:
\begin{align}
 G^{\mu,\nu}(t,x) := \sum_{j=0}^\infty{\Lambda_j^{\mu,\nu}(x)t^j}.\label{eq:GenFctExp}
\end{align}

\begin{fact}[{\cite[Theorem A]{HKMM09a}}]
$\Lambda_j^{\mu,\nu}(x)$ is real analytic on $\mathbb{R}_+$ and an $L^2$-solution of \eqref{eq:DiffEq}.
\end{fact}

We will now compute the constants $A_j^{\mu,\nu}(u)$ and $B_j^{\mu,\nu}(u)$ for $u=\Lambda_j^{\mu,\nu}$. Here we recall that $A_j^{\mu,\nu}$ and $B_j^{\mu,\nu}$ were defined in Theorem A and \eqref{eq:Asymptotics}.

From \cite[Theorem 4.2]{HKMM09a} we immediately obtain
\begin{align}
 B_j^{\mu,\nu}(\Lambda_j^{\mu,\nu}) &= \frac{\Gamma(j+\frac{\mu-|\nu|+2}{2})}{j!\Gamma(\frac{\mu+2}{2})\Gamma(\frac{\mu-|\nu|+2}{2})}\times \left\{\begin{array}{ll}\displaystyle2^{\nu-1}\Gamma\left(\frac{\nu}{2}\right)&\mbox{for $\nu>0$,}\\\displaystyle-1&\mbox{for $\nu=0$,}\\\displaystyle\frac{1}{2}\Gamma\left(-\frac{\nu}{2}\right)&\mbox{for $\nu=-1$.}\end{array}\right.\label{eq:BofLambda}
\end{align}

For $A_j^{\mu,\nu}$ we have the following lemma:

\begin{lemma}\label{lem:Ajmunu}
For any $j\in\mathbb{N}_0$ we have
\begin{align}
 A_j^{\mu,\nu}(\Lambda_j^{\mu,\nu}) &= (-1)^j\frac{2^{2(\mu+\nu)}\Gamma(j+\frac{\mu-\nu+2}{2})}{\pi\Gamma(j+\mu+1)}.\label{eq:AofLambda}
\end{align}
\end{lemma}

Putting \eqref{eq:BofLambda} and \eqref{eq:AofLambda} together proves Theorem B. In the remaining part of this section, we give a proof of Lemma \ref{lem:Ajmunu}.

\begin{proof}[\sc Proof of Lemma \ref{lem:Ajmunu}]
We put $\vartheta=\phi=0$ in \eqref{eq:MainIntFormula}. Using the special values
\begin{align*}
 \widetilde{J}_\alpha(0) &= \frac{1}{\Gamma(\alpha+1)}, & \widetilde{C}_n^\lambda(1) &= \frac{\Gamma(n+2\lambda)\Gamma(\lambda)}{\Gamma(n+1)\Gamma(2\lambda)},
\end{align*}
we obtain
\begin{align*}
 A_j^{\mu,\nu}(\Lambda_j^{\mu,\nu}) &= \frac{2^{\mu+\nu}j!\Gamma(j+\frac{\mu-\nu+2}{2})}{\pi\Gamma(j+\mu+1)\Gamma(j+\frac{\mu+\nu+2}{2})}\int_0^\infty{\Lambda_j^{\mu,\nu}(x)x^{\mu+\nu+1}\td x}.
\end{align*}
Together with the next lemma this finishes the proof.
\end{proof}

\begin{lemma}\label{lem:L1int}
For every $j\in\mathbb{N}_0$ we have $\Lambda_j^{\mu,\nu}\in L^1(\mathbb{R}_+,x^{\mu+\nu+1}\td x)$ and
\begin{align*}
 \int_0^\infty{\Lambda_j^{\mu,\nu}(x)x^{\mu+\nu+1}\td x} &= (-1)^j\frac{2^{\mu+\nu}\Gamma(j+\frac{\mu+\nu+2}{2})}{j!}.
\end{align*}
\end{lemma}

\begin{pf}
The fact that $\Lambda_j^{\mu,\nu}\in L^1(\mathbb{R}_+,x^{\mu+\nu+1}\td x)$ is derived from the asymptotic behaviour of $\Lambda_j^{\mu,\nu}(x)$ (see \cite[Theorem 4.2]{HKMM09a}). To calculate the integral we use the following integral formula which is valid for $\Re(\lambda+\alpha\pm\beta+1)>0$ and $b>a>0$ (see e.g. \cite[formula 6.576 (5)]{GR65}):
\begin{multline*}
 \int_0^\infty{I_\alpha(ax)K_\beta(bx)x^{\lambda}\td x} = \frac{a^\alpha\Gamma(\frac{\lambda+\alpha+\beta+1}{2})\Gamma(\frac{\lambda+\alpha-\beta+1}{2})}{2^{1-\lambda}b^{\lambda+\alpha+1}\Gamma(\alpha+1)}\times\\
 {_2F_1}\left(\frac{\lambda+\alpha+\beta+1}{2},\frac{\lambda+\alpha-\beta+1}{2};\alpha+1;\frac{a^2}{b^2}\right),
\end{multline*}
where ${_2F_1}(\alpha,\beta;\gamma;z)$ denotes the hypergeometric function. With \eqref{eq:GenFct} we obtain
\begin{align*}
 \int_0^\infty{G^{\mu,\nu}(t,x)x^{\mu+\nu+1}\td x} &= 2^{\mu+\nu}\Gamma\left(\frac{\mu+\nu+2}{2}\right)(1-t)^{\frac{\mu+\nu+2}{2}}\times\\
 & \ \ \ \ \ \ \ \ \ \ \ \ \ \ \ \ \ \ \ \ \ {_2F_1}\left(\frac{\mu+\nu+2}{2},\frac{\mu+2}{2};\frac{\mu+2}{2};t^2\right)\\
 &= 2^{\mu+\nu}\Gamma\left(\frac{\mu+\nu+2}{2}\right)(1+t)^{-\frac{\mu+\nu+2}{2}}\\
 &= \sum_{j=0}^\infty{\frac{2^{\mu+\nu}\Gamma(j+\frac{\mu+\nu+2}{2})}{j!}(-t)^j}.
\end{align*}
Then, in view of \eqref{eq:GenFctExp}, the claim follows by comparing coefficients of $t^j$.
\end{pf}

Hence, the proof of Theorem B is completed.

\section{Applications and special values}\label{sec:Applications}

We conclude this article with some applications of Theorem A and discuss on special values of the integral formula.

\subsection{The $L^2$-norm of $\Lambda_j^{\mu,\nu}$}

As a first application of Theorem A we can give a closed formula for the $L^2$-norms of the orthogonal basis $(\Lambda_j^{\mu,\nu}(x))_{j\in\mathbb{N}_0}$ in $L^2(\mathbb{R}_+,x^{\mu+\nu+1}\td x)$. The same result was obtained in \cite[Theorem B]{HKMM09a} by different methods.

\begin{corollary}
The $L^2$-norm of the eigenfunction $\Lambda_j^{\mu,\nu}$ is given by
 \begin{equation*}
  \|\Lambda_{2,j}^{\mu,\nu}\|_{L^2(\mathbb{R}_+,x^{\mu+\nu+1}\td x)}^2 = \frac{2^{\mu+\nu-1}\Gamma(j+\frac{\mu+\nu+2}{2})\Gamma(j+\frac{\mu-\nu+2}{2})}{j!(2j+\mu+1)\Gamma(j+\mu+1)}.
 \end{equation*}
\end{corollary}

\begin{pf}
Let $p=\mu+3$, $q=\nu+3$. We define functions $\phi_j\in L^2(C)$ in bipolar coordinates by
\begin{align*}
 \phi_j(r,\omega,\eta):=\Lambda_j^{\mu,\nu}(2r).
\end{align*}
Then it is immediate that
\begin{align*}
 \|\Lambda_j^{\mu,\nu}\|_{L^2(\mathbb{R}_+,x^{\mu+\nu+1}\td x)}^2 &= \frac{2^{\mu+\nu+3}}{\vol(S^{p-2})\vol(S^{q-2})}\|\phi_j\|_{L^2(C)}^2.
\end{align*}
Now, the intertwining operator $\mathcal{T}:L^2(C)\rightarrow\mathcal{H}$ is unitary up to a constant, namely (see \cite[Remark 2.2.2]{KM07b}):
\begin{align*}
 \|\mathcal{T}u\|_{\mathcal{H}}^2 &= \frac{1}{2}\|u\|_{L^2(C)}^2.
\end{align*}
Using formula \eqref{eq:RadialIntertwiner} for the intertwiner $\mathcal{T}$ and Theorem A one also obtains easily that
\begin{align*}
 \mathcal{T}\phi_j &= 2^{-\frac{3(\mu+\nu+2)}{2}}A_j^{\mu,\nu}(\Lambda_j^{\mu,\nu})\psi_j
\end{align*}
with $\psi_j$ as in \eqref{eq:KtypesConformal}. By \cite[Fact 2.1.1 (4)]{KM07b} the $G$-invariant norm on $\mathcal{H}$ is given by
\begin{align*}
 \|u\|^2_{\mathcal{H}} &= \left(j+\frac{p-2}{2}\right)\|u\|_{L^2(M)}^2, & \mbox{for }u\in V^j.
\end{align*}
By using the formula (see \cite[7.313 (2)]{GR65})
\begin{align*}
 \int_0^\pi{\left[\widetilde{C}_n^\lambda(\cos\vartheta)\right]^2\sin^{2\lambda}\vartheta\td\vartheta} &= \frac{\pi2^{1-2\lambda}\Gamma(n+2\lambda)}{n!(n+\lambda)},
\end{align*}
we get
\begin{align*}
 \|\psi_j\|_{L^2(M)}^2 &= \int_{S^{p-1}\times S^{q-1}}{\left|\widetilde{C}_j^{\frac{p-2}{2}}(v_0)\widetilde{C}_{j+\frac{p-q}{2}}^{\frac{q-2}{2}}(v_{p+q-1})\right|^2\td v}\\
 &= \vol(S^{p-2})\vol(S^{q-2})\int_0^\pi{\left[\widetilde{C}_j^{\frac{p-2}{2}}(\cos\vartheta)\right]^2\sin^{p-2}\vartheta\td\vartheta}\\
 &\ \ \ \ \ \ \ \ \ \ \ \ \ \ \ \ \ \ \ \ \ \ \ \ \ \ \ \ \ \ \ \ \ \ \ \times\int_0^\pi{\left[\widetilde{C}_{j+\frac{p-q}{2}}^{\frac{q-2}{2}}(\cos\phi)\right]^2\sin^{q-2}\phi\td\phi}\\
 &= \frac{\pi^2\Gamma(j+\mu+1)\Gamma(j+\frac{\mu+\nu+2}{2})}{j!2^{\mu+\nu}(j+\frac{\mu+1}{2})^2\Gamma(j+\frac{\mu-\nu+2}{2})}\vol(S^{p-2})\vol(S^{q-2}).
\end{align*}
Finally, putting all the steps together shows the claim.
\end{pf}

\subsection{The Poisson kernel of the Gegenbauer Polynomials}

As another application of Theorem A we get an integral formula for the generating function $G^{\mu,\nu}(t,x)$ of $L^2$-eigenfunctions, which is closely related to the Poisson kernel of the Gegenbauer polynomials. For this we put
\begin{multline*}
 I_{\mu,\nu}(t,\vartheta,\phi) := \int_0^\infty{\widetilde{I}_{\frac{\mu}{2}}\left(\frac{tx}{1-t}\right)\widetilde{K}_{\frac{\nu}{2}}\left(\frac{x}{1-t}\right)\widetilde{J}_{\frac{\mu}{2}}(ax)\widetilde{J}_{\frac{\nu}{2}}(bx)x^{\mu+\nu+1}\td x},
\end{multline*}
where $a:=\frac{\sin\vartheta}{\cos\vartheta+\cos\phi}$ and $b:=\frac{\sin\phi}{\cos\vartheta+\cos\phi}$.

\begin{corollary}\label{cor:IntFormulaGenFct}
For $\cos\vartheta+\cos\phi>0$ and $-1<t<1$ we have the following formula:
\begin{multline*}
 (1-t)^{-\frac{\mu+\nu+2}{2}}\left(\frac{\cos\vartheta+\cos\phi}{2}\right)^{-\frac{\mu+\nu+2}{2}}I_{\mu,\nu}(t,\vartheta,\phi)\\
 = \frac{2^{2(\mu+\nu)}}{\pi}\sum_{j=0}^\infty{(-1)^j\frac{\Gamma(j+\frac{\mu-\nu+2}{2})}{\Gamma(j+\mu+1)}\widetilde{C}_j^{\frac{\mu+1}{2}}(\cos\vartheta)\widetilde{C}_{j+\frac{\mu-\nu}{2}}^{\frac{\nu+1}{2}}(\cos\phi)t^j}.
\end{multline*}
\end{corollary}

For $\mu=\nu=2\lambda-1$ this yields a formula for the Poisson kernel of the Gegenbauer polynomials. Recall that the Poisson kernel $P_\lambda(t,\vartheta,\phi)$ of the Gegenbauer polynomials is defined by (see \cite[formula (6.4.5)]{AAR99})
\begin{align*}
 P_\lambda(t,\vartheta,\phi) := \sum_{n=0}^\infty{\frac{n!(n+\lambda)}{\Gamma(n+2\lambda)}\widetilde{C}_n^\lambda(\cos\vartheta)\widetilde{C}_n^\lambda(\cos\phi)t^n}.
\end{align*}
(For explicit formulas for the Poisson kernel of the Gegenbauer polynomials see e.g. \cite[formula (7.5.6)]{AAR99}.) Now, Corollary \ref{cor:IntFormulaGenFct} yields a new expression for $P_\lambda(t,\vartheta,\phi)$ ($2\lambda\in\mathbb{Z}$, $\lambda>0$):
\begin{multline*}
 P_\lambda(t,\vartheta,\phi) = \frac{\pi}{2^{8\lambda-4}}\left(\frac{\cos\vartheta+\cos\phi}{2}\right)^{-2\lambda}\\
 \times(\theta_t+\lambda)\left[(1+t)^{-2\lambda}I_{2\lambda-1,2\lambda-1}(-t,\vartheta,\phi)\right],
\end{multline*}
where $\theta_t=t\frac{\partial}{\partial t}$.

\subsection{The bottom layer}\label{sec:BottomCase}

As remarked in the introduction, for $j=0$ the $K$-Bessel function $u(x)=\widetilde{K}_{\frac{\nu}{2}}(x)$ is a solution of the differential equation \eqref{eq:DiffEq}. In this case the integral formula in Theorem A can be written as
\begin{multline*}
 \int_0^\infty{K_{\frac{\nu}{2}}(x)J_{\frac{\mu}{2}}\left(\frac{x\sin\vartheta}{\cos\vartheta+\cos\phi}\right)J_{\frac{\nu}{2}}\left(\frac{x\sin\phi}{\cos\vartheta+\cos\phi}\right)x^{\frac{\mu+2}{2}}\td x}\\
 = \frac{2^{\frac{\nu-2}{2}}\Gamma(\frac{\mu-\nu+2}{2})}{\sqrt{\pi}}(\cos\vartheta+\cos\phi)\sin^{\frac{\mu}{2}}\vartheta\sin^{\frac{\nu}{2}}\phi\ \widetilde{C}_{\frac{\mu-\nu}{2}}^{\frac{\nu+1}{2}}(\cos\phi).
\end{multline*}
This special case was already proved in \cite[Lemma 7.8.1]{KM07b}. Another expression for this integral can be found in \cite[formula 8.13 (14)]{EMOT54a}.

\addcontentsline{toc}{section}{References}

\providecommand{\bysame}{\leavevmode\hbox to3em{\hrulefill}\thinspace}
\providecommand{\MR}{\relax\ifhmode\unskip\space\fi MR }
\providecommand{\MRhref}[2]{%
  \href{http://www.ams.org/mathscinet-getitem?mr=#1}{#2}
}
\providecommand{\href}[2]{#2}

\end{document}